\documentclass[reqno]{amsart} 

\usepackage{amsmath}
\usepackage{amsfonts}
\usepackage{amssymb}
\usepackage{amsthm}
\usepackage[dvips]{epsfig}  
\usepackage{mathrsfs}
\usepackage[mathcal]{euscript}  
\usepackage{natbib}  
 
\numberwithin{equation}{section}   
   
\theoremstyle{plain}   
\newtheorem{thm}{Theorem}[section]   
\newtheorem{lem}[thm]{Lemma}

\theoremstyle{definition}   

\newtheorem{exmp}[thm]{Example}
\theoremstyle{remark}   
\newtheorem{rem}[thm]{Remark}
\newcommand{\eqdis}{\stackrel{\lower0.2ex\hbox{$\scriptscriptstyle   
                    \mathrm{d}$}}{=}}   
\newcommand{\vague}{\stackrel{\lower0.2ex\hbox{$\scriptscriptstyle   
                    \it{v} $}}{\rightarrow}}   
\newcommand{\weak}{\stackrel{\lower0.2ex\hbox{$\scriptscriptstyle   
                    \it{w} $}}{\rightarrow}}   
\newcommand{\what}{\stackrel{\lower0.2ex\hbox{$\scriptscriptstyle   
                    \it{\hat{w}} $}}{\rightarrow}}   
\newcommand{\distr}{\stackrel{\lower0.2ex\hbox{$\scriptscriptstyle   
                    \it{d} $}}{\rightarrow}}   
\newcommand{\inprob}{\stackrel{\lower0.2ex\hbox{$\scriptscriptstyle   
      \Prob$}}{\rightarrow}}   
\newcommand{\Prob}{{P}}
\newcommand{\Var}{\operatorname{Var}}   
   
\newcommand{\R}{\mathbf{R}}   
\newcommand{\vep}{\varepsilon}

\newtheorem{alg}{Algorithm}

 \makeatletter
 \AtBeginDocument{%
  \def\@serieslogo{%
  \vbox to\headheight{%
  \parindent\z@ \fontsize{6}{7\p@}\selectfont
  September 11, 2009\endgraf
  \vss}}}
 \makeatother

\begin{document}   
 
\title[Importance sampling for random walks with heavy tails]{On
  importance sampling with mixtures for random walks 
  with heavy tails}

\author[H.~Hult]{Henrik Hult}  
\address[H.~Hult]{Department of Mathematics, KTH, 100 44 Stockholm, Sweden} 
\email{hult@kth.se}  
\author[J.~Svensson]{Jens Svensson}  
\address[J.~Svensson]{Department of Mathematics, KTH, 100 44 Stockholm, Sweden} 
\email{jenssve@kth.se}

\copyrightinfo{}{The authors}

\begin{abstract}   
  Importance sampling algorithms for heavy-tailed random
  walks are considered. Using a specification with  
  algorithms based on mixtures of the original distribution
  with some other distribution, sufficient conditions 
  for obtaining bounded relative error are presented.
  It is proved that mixture algorithms of this kind 
  can achieve asymptotically optimal relative error.
  Some examples of mixture algorithms are presented, including
  mixture algorithms using a scaling of the original 
  distribution, and the bounds of the relative errors 
  are calculated. The algorithms are evaluated numerically in a simple setting.
\end{abstract}   
   

\maketitle   

\section{Introduction}
Tail probabilities appear naturally in many applications of probability theory, 
and often analytical evaluation is not possible. For many applications, 
Monte Carlo simulation can be an effective alternative. For rare events, 
however, standard Monte Carlo simulation is very computationally inefficient, 
and some form of variance reduction method is necessary. One such alternative 
that has been extensively applied to both light- and heavy-tailed
distributions is importance sampling. 
In this paper we focus on importance sampling algorithms for
computing the probability 
\begin{align*}
  p_b = P(S_n > b),
\end{align*}
of a high threshold $b$, for a random walk  $S_n = X_1 + \dots +
X_n$. The random variables
$X_1,\dots, X_n$ are independent and identically distributed with 
distribution function $F$ and density $f$. It
is assumed that the   
right tail of $f$ is regularly varying at $\infty$; more precisely there exists an 
$\alpha > 0$ such that, for each $x > 0$,
\begin{align*}
  \lim_{u \to \infty}\frac{f(ux)}{f(u)} = x^{-\alpha-1}.
\end{align*}
Then it is well known that $f$ has the representation $f(x) =
x^{-\alpha-1}L(x)$, $x > 0$, where $L$ is slowly varying. The joint
distribution of $(X_1,\dots,X_n)$ is denoted $\mu_n$.

Consider first a computation of $p_b$ using standard Monte Carlo. Then $N$
independent samples $(X_1^1, \dots, X_n^1), \dots, (X_1^N, \dots,
X_n^N)$ are generated from $\mu_n$ and $p_b$ is estimated using the
sample frequency
\begin{align*}
  \hat p_b^{MC} = \frac{1}{N} \sum_{i=1}^N I\{S_n^i > b\},
\end{align*}
where $S_n^i = X_1^i+\dots+X_n^i$. For large $b$, the event
$\{S_n>b\}$ is rare and few of the indicator variables $I\{S_n^i >
b\}$ will be $1$. This leads to rather inefficient
estimation. To see this, consider for instance the standard deviation of $\hat
p_b^{MC}$. An elementary calculation shows
\begin{align*}
  \text{Stdev}(\hat p_b) = \frac{1}{\sqrt{N}}\sqrt{p_b(1-p_b)}.
\end{align*}
When $p_b$ is small this is roughly $\sqrt{p_b/N}$. Hence, it would
require $N \approx 1/p_b$ samples to have the standard deviation of
size comparable to the quantity $p_b$ we are estimating. When $p_b$ is
small this can be very large. 

Importance sampling provides a way to possibly reduce the
computational cost without sacrificing precision, or equivalently to
improve precision without increasing the computational cost. 
The basic idea of importance sampling to generate samples $(X_1^1,
\dots, X_n^1), \dots, (X_1^N, \dots, 
X_n^N)$ independently from a sampling measure $\nu_n^b$ instead of
$\mu_n$. It is assumed that $\mu_n$ is absolutely continuous with
respect to $\nu_n^b$, written $\mu_n \ll \nu_n^b$ so that the Radon-Nikodym 
derivative $\frac{d\mu_n}{d\nu_n^b}$ exists. An unbiased estimate
of $p_b$ is constructed as 
\begin{align*}
  \hat p_b = \frac{1}{N} \sum_{i=1}^N
  \frac{d\mu_n}{d\nu_n^b}(X_1,\dots,X_n)I\{S_n > b\}.
\end{align*}
The goal is to choose $\nu_n^b$ so more samples 
are drawn from regions that are  ``important'' to the event
$\{S_n>b\}$. Then the event becomes less rare under $\nu_n^b$, which reduces 
variance. However, $\nu_n^b$ must be chosen carefully so that the
Radon-Nikodym weights $\frac{d\mu_n}{d\nu_n^b}(X_1,\dots,X_n)$ do not cause
variance to increase. A relevant quantity for deciding if a
sampling measure $\nu_n^b$ is appropriate or not is 
the relative error
\begin{align*}
  RE(\hat p_b) = \frac{\sqrt{\Var(\hat p_b)}}{E \hat p_b} =
  \frac{1}{\sqrt{N}}\sqrt{\frac{E \hat p_b^2 - p_b^2}{p_b^2}} =
  \frac{1}{\sqrt{N}}\sqrt{\frac{E \hat p_b^2}{p_b^2} - 1}. 
\end{align*}
By Jensen's inequality we always have $E \hat p_b^2 \geq p_b^2$. 

To quantify the efficiency of the sampling measure it is convenient to study
the asymptotics of the relative error as $b \to \infty$. This amounts
to studying the asymptotics of normalized second moment $\lim_{b \to
  \infty} E\hat p_b^2/p_b^2$. We say that
a sampling distribution $\nu_n^b$ has 
{\it logarithmically efficient relative error} if, for some $\vep > 0$,
\begin{align*}
  \limsup_{b \to \infty} \frac{E \hat p_b^2}{p_b^{2-\vep}} < \infty,
\end{align*}
it has {\it bounded relative error}  
if
\begin{align*}
  \limsup_{b \to \infty} \frac{E \hat p_b^2}{p_b^{2}} < \infty,
\end{align*}
and {\it asymptotically optimal relative error} if
\begin{align*}
  \limsup_{b \to \infty} \frac{E \hat p_b^2}{p_b^2} =1. 
\end{align*}

A number of different algorithms have been proposed to simulate tail
probabilities of heavy-tailed random walks. \cite{AB97} study
the class of subexponential distributions,  i.e.~distributions for which 
\begin{align*}
  \lim_{b\to\infty}\frac{P(S_n>b)}{nP(X_1>b)}=1,
\end{align*}
and use that as $b\to\infty$, all the variance of the sum comes from
the largest summand. By removing the largest term $X_{(n)}$ in each
sample and calculating the probability using the remaining $n-1$
terms, they obtain a logarithmically efficient conditional Monte Carlo
estimator in the sub-class of distributions with regularly varying
tails. Here $X_{(1)}<X_{(2)}<\ldots<X_{(n)}$ is the ordered sample. 
Specifically, with the sample $X_1, \ldots, X_n$, the estimator is 
\begin{align*}
  P(S_n>b|X_{(1)}, \ldots , X_{(n-1)})=\frac{ \overline{F} (X_{(n-1)}\vee (b-S_{(n-1)}))}{\overline{F}(X_{(n-1)})},
\end{align*}
where $S_{(n-1)}=X_{(1)}+\ldots X_{(n-1)}$ is the sum of the $n-1$ largest of the 
$n$ terms in the sample and $a\vee b=\max(a,b)$.

\cite{AK06} propose a similar idea, using the conditioning
\begin{align*}
  nP(S_n>b, M_n=X_n|X_{1}, \ldots , X_{n-1})= \overline{F} (M_{n-1}\vee (b-S_{(n-1)})),
  \end{align*}
where $M_n=\max(X_1,\ldots, X_n)$ to obtain an estimator with bounded relative 
error for distributions with regularly varying tails.

\cite{JS02} introduce an importance sampling algorithm with 
similar structure to  exponential twisting in the light-tailed case. 
Their so-called hazard rate twisting of the original distribution is given by 
\begin{align}
  dF_{\theta}(x)=
  \frac{e^{\theta\Lambda(x)}dF(x)}{\int_{0}^{\infty}e^{\theta\Lambda(x)}dF(x)},  
\end{align}
 where $0<\theta<1$ and $\Lambda(x)=-\log \overline{F}(x)$ is the hazard rate. 
For distributions with regularly varying tails, this is equivalent to changing 
the tail index of the distribution.
 
In the case of importance sampling, \cite{BJZ07} show that to obtain 
efficient sampling distributions for heavy-tailed random walks, one
must consider  
state-dependent changes of measure. Simply changing the
parameters in the original distribution cannot lead to an 
estimator with bounded relative error. 

The first algorithm of this type, for heavy-tailed random walks,  was
proposed by \cite{DLW07}. There the large values are sampled
from the conditional distribution where one has to condition on
exceeding a level just below the remaining distance to $b$. The
authors prove that their proposed algorithm has close to
asymptotically optimal relative error. 
\cite{BLi08} present a state-dependent algorithm that uses Markov 
chain description of the random walk under the sampling measure
to obtain bounded relative error for the class of subexponential
distributions.
\cite{BLiu08} construct a mixture algorithm with bounded relative error
for the large deviation probability $P(S_n>b)$ where $b>b_0n^{1/2+\epsilon}$.

In this paper we take a more general look at mixture algorithms of the
same type as \cite{DLW07}.
The underlying idea is to construct a dynamic change of measure such that the
trajectories of $X_1, \dots, X_n$ leading to $S_n > b$ is similar to the
most likely trajectories conditional on $S_n > b$. In the heavy-tailed
case, the most likely trajectories are such that one of the $X_i$'s is 
large and the others are ``average''. Mixtures arise quite naturally as 
sampling distributions for producing such trajectories; with some probability 
$p_i$ sample from the original density $f$ and with
probability $q_i = 1-p_i$ sample from a density where it is likely to
get a large value.  
We provide sufficient conditions for bounded relative error and 
provide a couple of new examples that are very easy
to implement. We also show that, with some additional work, one can construct
mixture algorithms with asymoptotically optimal relative error. 

The paper is organized as follows. In Section 2 we
present a general importance sampling algorithm based on mixtures and
provide several examples. In Section 3 we provide
sufficient condition for the mixture algorithm to have bounded
relative error. In Section 4 we provide detailed analysis
of specific mixture algorithms. The concluding Section 5 
provides a proof that it is possible to obtain asymptotically optimal 
relative error. 

\section{Dynamic mixture algorithms}\label{sec:mixalg}

In this section we describe a general importance sampling algorithm
based on mixtures, called the dynamic mixture algorithm, and provide
several examples. 

The dynamic mixture algorithm for computing $p_b = P(S_n > b)$
proceeds as follows. Each replication of $(X_1, \dots, X_n)$ is
generated dynamically and the distribution for sampling $X_i$ depend
on the current state $S_{i-1} = X_1 + \dots + X_{i-1}$ of the random
walk. At the $i$th step it may be that $S_{i-1}$ already exceeds the
threshold $b$. Then $X_i$ is sampled from the 
original density $f$. Otherwise, if $S_{i-1} \leq b$, a biased coin is
tossed with probability $p_i$ for ``heads'' and $q_i = 1-p_i$ for
``tails''. If it comes up ``heads'' $X_i$ is generated from the usual
density $f$, but if it comes up ``tails'', $X_i$ is generated from another
density $g_i(x \mid S_{i-1})$.  The density $g_i(x \mid S_{i-1})$
depends on the current generation $i$ of the algorithm 
and on the current position 
$S_{i-1}$. The idea is to choose $g_i(x \mid S_{i-1})$ s.t.\ sampling 
from $g_i(x \mid S_{i-1})$ is likely to result in a large variable. 
 However, $g_i(x \mid S_{i-1})$ must be chosen with some care to
control the Radon-Nikodym weights $\frac{d\mu_n}{d\nu_n^b}(X_1, \dots,
X_n)$ and thereby the relative error. In the last
generation, if $S_{n-1} \leq b$, $X_n$ is sampled from a density
$g_n(x \mid S_{n-1})$ and if $S_{n-1} > b$ it is sampled from the
original $f$. In contrast to the previous steps $g_n$ is not 
necessarily of mixture type. The reason is that it may be advantageous to
make sure $X_n > b-S_{n-1}$ in the last step to get $S_n > b$.

A precise description of the dynamic mixture algorithm is presented
next. 
\begin{alg}\label{alg:gm}
Consider step $i=1,\dots,n$, where $S_{i-1} = s_{i-1}$. Then
$X_i$ is sampled as follows.
\begin{itemize}
\item If $s_{i-1} > b$, $X_i$ is sampled from the original
  density $f$,
\item if $s_{i-1} \leq b$, $X_i$ is sampled from 
  \begin{align*}
    &p_if(\cdot) + q_i g_i( \cdot \mid s_{i-1}), \quad \text{ for }1
    \leq i \leq n-1,\\ 
    &g_n(\cdot \mid s_{n-1}), \quad \text{ for } i = n.
\end{align*}
Here $p_i + q_i = 1$ and $p_i \in (0,1)$. 
\end{itemize}
\end{alg}

Explicit examples of the dynamic mixture algorithm are obtained by
specifying $g_i$ and $p_i$.
\begin{exmp}[Conditional mixture, c.f.\ \cite{DLW07}] \label{exmp:cm}
  The algorithm proposed by \cite{DLW07} takes $g_i$ to be  a
  conditional distribution. For $i = 1,\dots, n-1$ the large values
  are sampled conditional on being at least $a$ times the remaining
  distance to $b$, where $a \in (0,1)$. It is important that
  $a<1$. In the last step samples are generated conditional on exceeding $b$. 
  More precisely, 
  \begin{align*}
    g_i(x \mid s) &= \frac{f(x)I\{x > a(b-s)\}}{\overline{F}(a(b-s))},
    \quad 1 \leq i \leq n-1,\\
    g_n(x \mid s) &= \frac{f(x)I\{x > b-s\}}{\overline{F}(b-s)}.
  \end{align*}
  In their paper the authors assume that $f = 0$ on $(-\infty,
  0)$. That is, all the $X_i$'s are non-negative. This is not an important 
  restriction and we do not impose it here. 
\end{exmp}

A practical limitation of the conditional mixture algorithm is that
some distributions do not allow direct sampling from the conditional  
distribution. If the distribution function $F$ and its inverse
$F^{\leftarrow}$ are available, the inversion method suggest sampling
$X$ conditional on $X > c$ by taking $U$ to be uniform on $(0,1)$ and
set $X = F^{\leftarrow}(1-U\overline{F}(x))$, see e.g.\ \cite{AG07}. In
other cases it might 
be necessary to use an acceptance-rejection method, but this may be
time consuming. 

A simple alternative to the conditional mixture is
to sample the large variables from a generalized Pareto distribution
(GPD) instead. The intuition is that the GPD approximates the
conditional distribution well.

\begin{exmp}[Generalized Pareto mixture] \label{exmp:gpd}
  The GPD mixture algorithm takes $g_i$ to be  a generalized Pareto
  distribution. As in the previous algorithm, for $i = 1,\dots, n-1$,
  the large values are sampled conditional on being at least $a$ times
  the remaining distance to $b$, where $a \in (0,1)$. The last step is
  slightly different. If the remaining distance is 
  large, the last step is taken from a GPD, otherwise it is taken from
  the original density. This is because, if $S_{n-1} \leq b$, but
  close to $b$, the GPD is not necessarily a good approximation of the
  conditional distribution. To be precise,
  \begin{align*}
    g_i(x \mid s) &= \alpha[a(b-s)]^\alpha x^{-\alpha-1} I\{x > a(b-s)\},
    \quad 1 \leq i \leq n-1.\\
    g_n(x \mid s) &= \alpha
    (b-s)^\alpha x^{-\alpha-1} I\{x > b-s\}I\{s \leq b-b(1-a)^{n-1}\}\\ 
    & \quad + f(x)I\{s > b-b(1-a)^{n-1}\}.
  \end{align*}
\end{exmp}

A different way to sample the large variables is to sample from
the original density and then scale the outcome by simply multiplying
with a large number $\lambda b$. We 
call this a scaling mixture algorithm.
\begin{exmp}[Scaling mixtures]\label{exmp:scalemix}
  The scaling mixture algorithm has, with
  $\lambda > 0$, 
  \begin{align*}
    g_i(x\mid s) &= (\lambda b)^{-1}f(x/\lambda b)I\{x > 0\} + f(x)I\{x
    \leq 0\},\quad i=1,\dots,n-1,\\
    g_n(x \mid s) &=  (\lambda b)^{-1}f(x/\lambda b)I\{x > 0, s \leq
    b-b(1-a)^{n-1}\}\\  
    & \quad + f(x)I\{x \leq 0 \text{ or } s > b-b(1-a)^{n-1}\}.
  \end{align*}
  To simplify the analysis we will, in the context of
  scaling mixtures,  always assume that the orginal density $f$ is
  strictly positive on $(0,\infty)$. If this is not satisfied the 
  situation is more involved because there may be large $x> 0$ such
  that $f(x) > 0$ but $f(x/\lambda b) = 0$. Then such large
  $x$-value cannot be obtained by sampling a small number from $f$ and
  scale by $\lambda b$. This may cause the Radon-Nikodym weights to be
  relatively large, which increase the variance. 
  
  There are several variations of the scaling algorithm. For instance,
  one may scale with something proportional to the remaining distance to
  $b$, instead of something proportional to $b$ as described above. 
  Some variations of the scaling algorithm will be treated in
  more detail in Section \ref{sec:scalingmix}.
\end{exmp}

\section{Asymptotic analysis of the normalized second moment}\label{sec:asm}

The efficiency criteria presented in the introduction are all based on
the asymptotic properties of the normalized second moment $E \hat
p_b^2/p_b^2$.  We are following the weak
convergence approach initiated by \cite{DLW07} to study its asymptotics. 
By the subexponential property, $p_b^2\sim n^2 \overline{F}(b)^2$, where $a_b \sim c_b$ denotes
$\lim_{b\to\infty} a_b/c_b=1$, the normalized second moment can be written as
\begin{align}\label{eq:nsm}
  \frac{E\hat p_b^2}{p_b^2} \sim 
  \frac{1}{n^2\overline{F}(b)^2}\int\limits_{s_n>b}\frac{d\mu_n}{d\nu_n^b}(y)
  \mu_n(dy)= 
  \frac{1}{n^2}\int\limits_{s_n>1}\frac{1}{\overline{F}(b)}\frac{d\mu_n}{d\nu_n^b}(by)m_b(dy),  
\end{align}
where the measure $m_b = \mu_n(b(\,\cdot\cap \{s_n>1\}))/\overline{F}(b)$.
To calculate the limit of this integral we will use the weak
convergence of the measure $m_b$ to a measure $m$ and uniform convergence of
an upper bound
$R^*_b(y) \geq 
\frac{1}{\overline{F}(b)}\frac{d\mu_n}{d\nu_n^b}(by) =:R_b(y)$ to a bounded  
continuous function $R(y)$. Then we (?) establish the convergence
\begin{align*}
  \limsup_{b \to \infty} \int\limits_{\{s_n>1\}} R_b dm_b \leq   \lim_{b
    \to \infty} \int\limits_{\{s_n>1\}} R_b^* dm_b  =   \int\limits_{\{s_n>1\}} R dm. 
\end{align*}
To do this it is convenient if the normalized Radon-Nikodym
derivative $R_b(y)$ is bounded. This criteria is certainly stronger
than necessary but appears to be desirable. It implies the the
normalized $q$-moment is asymptotically bounded for any $q \in
(1,\infty)$. Indeed, if $R_b^*$ is
bounded and $R^* \to R$ uniformly, then for any $q \in (1,\infty)$
\begin{align*}
  \limsup_{b \to \infty}\frac{E\hat p_b^q}{p_b^q} = \limsup_{b \to \infty}
  \frac{1}{n^q}\int\limits_{s_n>1}\Big(\frac{1}{\overline{F}(b)}\frac{d\mu_n}{d\nu_n^b}(by)\Big)^{q-1}m_b(dy)     
  \leq \frac{1}{n^q}\int R^{q-1} dm < \infty.
\end{align*}

Next we provide sufficient
conditions for $R_b$ to be bounded. 

\begin{lem}\label{lem:RNbddgen}
  Consider Algorithm \ref{alg:gm} with $p_i > 0$ for $1 \leq i
  \leq n-1$. Suppose there exists $a \in (0,1)$ such that
  \begin{align}
    \liminf_{b \to \infty}\inf_{\scriptsize \begin{array}{l} s \leq
        1-(1-a)^{i-1}\\ y > a(1-s)\end{array}}\frac{g_i(by\mid b s)}{f(by)} 
    \overline{F}(b) &> 0, \quad  1 \leq i \leq n, \label{eq:main}\\
    \limsup_{b \to \infty}\sup_{\scriptsize \begin{array}{l} s \leq 1\\ y >
        1-s \end{array}}\frac{f(by)}{g_n(by\mid b s)} &< \infty. \label{eq:sec}
  \end{align}
  Then the scaled Radon-Nikodym derivative $R_b(y) =
  \frac{1}{\overline{F}(b)}\frac{d\mu}{d\nu_n^b}(by)$ is bounded on
  $\{y_1+\dots + y_n > 1\}$.    
\end{lem}

\begin{proof}
  Let $s_n = y_1 + \dots + y_n$. On $\{s_n >
  1\}$ it must hold that $y_i > a(1-s_{i-1})$ for some $i =
  1,\dots,n$. Otherwise $s_i \leq 1-(1-a)^{i} < 1$ for each $i$. 
  
  Take $y \in \{y\in \mathbb{R}^n: s_n > 1\}$ and let $i = \min\{j: y_j >
  a(1-s_{j-1})\}$. Note that for this $i$
  \begin{align*}
    y_i > a(1-s_{i-1}) \geq a(1-a)^{i-1} \geq a(1-a)^n =: a_n > 0. 
  \end{align*}
  For any $y_j$, $j \notin \{i,n\}$, 
  \begin{align*}
     \frac{f(by_j)}{p_j f(by_j) + q_jg_j(by_j\mid b s_{j-1})} \leq \frac{1}{p_j}.
  \end{align*}
  It follows that, for $1 \leq i \leq n-1$,
  \begin{align}\label{eq:threetermsgen}
    \frac{1}{\overline{F}(b)}\frac{d\mu}{d\nu_n^b}(by) &\leq
    \frac{1}{\overline{F}(b)}\frac{f(by_i)I\{y_i >
      a(1-s_{i-1})\}}{p_i f(by_i) + q_ig_i(by_i \mid bs_{i-1})} \times 
    \prod_{j \notin \{i,n\}} \frac{1}{p_j} \nonumber \\
    & \quad \times \Big(\frac{f(by_n)I\{y_n>
      1-s_{n-1}\}}{g_n(by_n \mid bs_{n-1})}I\{s_{n-1} \leq 1\} +
    I\{s_{n-1} > 1\}\Big).
  \end{align}
  The first term can be written as
  \begin{align*}
    \frac{1}{\overline{F}(b)}\frac{f(by_i)I\{y_i >
      a(1-s_{i-1})\}}{p_if(by_i) + q_ig_i(by_i\mid b s_{i-1})} & = 
    \frac{I\{y_i > a(1-s_{i-1})\}}{p_i\overline{F}(b)
      + q_i\frac{g_i(by_i \mid b s_{i-1})}{f(b y_i)} \overline{F}(b)}.
  \end{align*}
  By \eqref{eq:main} this term in \eqref{eq:threetermsgen} is
bounded. The second term is bounded because $p_j > 0$ by assumption. 
  The last term is bounded by \eqref{eq:sec}. 

Similarly for $i = n$,  
\begin{align*}
  \frac{1}{\overline{F}(b)}\frac{d\mu}{d\nu_n^b}(by) &\leq
  \frac{f(by_n)I\{y_n > a(1-s_{n-1})\}}{g_n(by_n \mid bs_{n-1})\overline{F}(b)}
  \prod_{j=1}^{n-1} \frac{1}{p_j},
\end{align*}
which is bounded by \eqref{eq:main}.
\end{proof}

Next we present the main result. It provides sufficient conditions
for the mixture algorithms to have bounded relative error. This is
obtained by showing that the normalized second moment remains
bounded. 
\begin{thm}\label{thm:main}
Suppose \eqref{eq:main} and \eqref{eq:sec} hold for $a\in(0,1)$. Suppose in 
addition that there exist continuous functions $h_i: \R^n \to
[0,\infty)$ such that
  \begin{align}\label{eq:thm}
    \frac{f(by_i)}{g_i(by_i \mid bs_{i-1}) \overline{F}(b)} \to h_i(y_i
    \mid s_{i-1}), 
  \end{align}
uniformly on $\{y \in \R^n: s_{i-1} \leq 1-(1-a)^{i-1}, y_i > a(1-s)\}$. 
Then,  
\begin{align*}
  \lim_{b\to\infty} \frac{E \hat p_b^2}{p_b^2} \leq \frac{1}{n^2}
  \sum_{i=1}^n \prod_{j=1}^{i-1}\frac{1}{p_j} \frac{1}{q_i}
  \int_1^\infty  h_i(y_i \mid 0) \alpha y_i^{-\alpha-1}dy_i,
\end{align*}
with the convention that $q_n = 1$. 
\end{thm}

\begin{proof}
  First rewrite the normalized second moment as in \eqref{eq:nsm}:
  \begin{align*}
    \frac{E \hat p_b^2}{p_b^2} \sim  
    \frac{1}{n^2}\int\limits_{s_n>1}
    \frac{1}{\overline{F}(b)}\frac{d\mu_n}{d\nu_n^b}(by)m_b(dy) =
    \frac{1}{n^2}\int 
    R_b(y) m_b(dy).   
  \end{align*}
  By regular variation of $f$ and independence of $X_1, \dots, X_n$
  the joint distribution $\mu_n$ is multivariate regularly varying. In
  particular the weak convergence $m_b \weak m$ holds, where $m$ has the  
  representation
  \begin{align}\label{eq:m}
    m(A) = \sum_{i=1}^n \int_A I\{y \in \R^n: y_i > 1,
    y_j = 0, j \neq i\} \alpha y_i^{-\alpha-1}dy_i.
  \end{align}
  This is well known, see e.g.\ \cite{R87}, Section 5.5. A
  proof is also given by \cite{DLW07}. 
  We see that the measure $m$  puts all its 
  mass on the coordinate axes. That is, on trajectories where one jump
  is large and the rest are zero. 

  The next step is to decompose the integral as 
  \begin{align}\label{eq:intdec}
    \int R_b dm_b =  \int_A
    R_b dm_b  + \int_{A^c} R_b dm_b,
  \end{align}
  where $A = \cup_{i=1}^n A_i$ is a finite union and the $A_i$'s have
  disjoint closures. We will find $A_i$ such that the second integral
  converges to $0$ and determine an upper bound $R^*_b \geq R_b$
  on $A_i$.  

  Define the sets $A_i$ to be
  \begin{align*}
    A_i=\{y\in \mathbb{R}^n: y_j\leq a(1-s_{j-1})
    \textrm{ for } 1 \leq j \leq i-1, y_i>1-s_{i-1}, & \\  \text{and }
    s_k > 1, k = i+1,\dots,n &\}.
  \end{align*}  
  Note that the $A_i$'s have disjoint closure and $m(\partial A_i) =
  0$. In particular $m_b(A_i) \to m(A_i)$ for each $i = 1,\dots, n$.
  Moreover, $m(\cap_{i=1}^n A_i^c) = 0$. Indeed, 
  \begin{align*}
    m(\cap_{i=1}^n A_i^c) &= m(\{s_n > 1\} \setminus \cup_i A_i) 
    = m\{s_n > 1\} - \sum_{i=1}^n m(A_i) \\
    &= \sum_{i=1}^n \int_{\{s_n > 1\}}I{\{y \in \R^n: y_i > 1, y_j =
      0, j \neq i\}}\alpha y_i^{-\alpha-1}dy_i\\ & \quad  - 
    \sum_{i=1}^n \int_{A_i} I{\{y \in \R^n: y_i > 1, y_j = 0, j \neq i\}}
    \alpha y_i^{-\alpha-1}dy_i = 0.
  \end{align*}
  By Lemma \ref{lem:RNbddgen}, $R_b$ is bounded and since $m_b(A^c) \to
  m(A^c) = 0$, the second integral in \eqref{eq:intdec} converges to $0$. 
  For the first integral we construct a function $R_b^*$ that
  dominates $R_b$ on $A$ and a continuous function $R$ such that
  $R_b^* \to R$ uniformly on $A$. Then it follows from weak
  convergence that
  \begin{align*}
    \limsup_{b\to \infty} \int_A R_b dm_b \leq  \lim_{b\to \infty} \int_A
    R_b^* dm_b  = \int_A R dm < \infty.  
  \end{align*}
  For $y \in A_i$,
  \begin{align*}
    R_b(y) &\leq \frac{1}{\overline{F}(b)}\frac{f(by_i)I\{y_i >
      1-s_{i-1}\}}{p_i f(by_i) + q_ig_i(by_i \mid bs_{i-1})} 
    \prod_{j=1}^{i-1} \frac{1}{p_j} =: R_b^*(y).
  \end{align*}
  To see this, construct a bound as in \eqref{eq:threetermsgen}
  and notice that on $A_i$, $s_k > 1$ for each $k \geq i$. Then the
  contribution to the Radon-Nikodym weights from $y_k$, $k > i$ is equal to
  $1$. By assumption \eqref{eq:thm}
  \begin{align*}
    \frac{g_i(by_{i} \mid bs_{i-1})}{f(by_i)} \overline{F}(b) \to \frac{1}{h_i(y_i
    \mid s_{i-1})}, 
  \end{align*}
uniformly on $A_i$. For $y \in A_i$ define $R(y) = 
  h_i(y_i\mid s_{i-1})\prod_{j=1}^{i-1}\frac{1}{p_j} \frac{1}{q_i}$. Then
$R_b^* \to R$ uniformly on $A$.    
With the representation \eqref{eq:m} of the limiting measure $m$, 
the upper bound for the normalized second moment can now be calculated as 
\begin{align*}
  \frac{1}{n^2}\int_A R dm =  \frac{1}{n^2}\sum_{i=1}^n
  \prod_{j=1}^{i-1}\frac{1}{p_j} \frac{1}{q_i} 
  \int_1^\infty h_i(y_i \mid 0)\alpha y_i^{-\alpha-1} dy_i.
\end{align*}
\end{proof}

\section{Examples}\label{sec:exmp}
In this section we provide a detailed analysis of the algorithms
presented in Section \ref{sec:mixalg}. In particular
we verify the conditions of Lemma \ref{lem:RNbddgen} and Theorem
\ref{thm:main} for these algorithms.

\subsection{The conditional mixture algorithm}  
Recall from Example \ref{exmp:cm} that the conditional mixture
algorithm  has, with $a \in (0,1)$,
\begin{align*}
  g_i(x \mid s) &= \frac{f(x)I\{x > a(b-s)\}}{\overline{F}(a(b-s))},
  \quad 1 \leq i \leq n-1,\\
  g_n(x \mid s) &= \frac{f(x)I\{x > b-s\}}{\overline{F}(b-s)}.
\end{align*}
Then, for $i = 1,\dots, n-1$, the uniform convergence 
$\overline{F}(bx)/\overline{F}(b) \to x^{-\alpha}$,  for $x > x_0 >
0$, implies
\begin{align*}
  \frac{g_i(bx \mid bs)}{f(bx)}\overline{F}(b) &=
  \frac{\overline{F}(b)}{\overline{F}(ba(1-s))}I\{x > a(1-s)\} \to
  a^{\alpha}(1-s)^{\alpha}I\{x > a(1-s)\},
\end{align*}
uniformly for $s \leq 1-(1-a)^{i-1}$, $x > a(1-s)$. Similarly, 
\begin{align*}
  \frac{g_n(bx \mid bs)}{f(bx)}\overline{F}(b) &=
  \frac{\overline{F}(b)}{\overline{F}(b(1-s))}I\{x > 1-s\} \to
  (1-s)^{\alpha}I\{x > 1-s\},
\end{align*}
uniformly on $s \leq 1-(1-a)^{n-1}$, $x > 1-s$,
and
\begin{align*}
  \frac{f(bx)}{g_n(bx \mid bs)} = \overline{F}(b(1-s)) \leq 1,
\end{align*}
on $s \leq 1$.
It follows that both \eqref{eq:main} and \eqref{eq:sec} are satisified
and hence the normalized Radon-Nikodym derivative is bounded. 

By the above calculation \eqref{eq:thm} holds with $h_i(y \mid s) =
a^{-\alpha} (1-s)^{-\alpha}$, $1 \leq i \leq n-1$ and $h_n(y \mid s) =
(1-s)^{-\alpha}$. It follows from Theorem \ref{thm:main} that 
\begin{align}\label{eq:condmix}
   \lim_{b\to \infty} \frac{E \hat p_b^2}{p_b^2} \leq \frac{1}{n^2}
   \int R(y)dm =  \frac{1}{n^2} \Big(\sum_{i=1}^{n-1} 
   \frac{a^{-\alpha}}{q_i}\prod_{j=1}^{i-1}\frac{1}{p_j}  +
   \prod_{j=1}^{n-1}\frac{1}{p_j}\Big). 
\end{align}
The right hand side is minimized at 
\begin{align}\label{eq:p}
  p_i = \frac{(n-i-1)a^{-\alpha/2}+1}{(n-i)a^{-\alpha/2}+1},\quad q_i
  = 1-p_i,
\end{align}
with minimum $n^{-2}[(n-1)a^{-\alpha/2}+1]^2$, and it is possible
to show that the limit is equal to the right hand side of \eqref{eq:condmix}, 
see \cite{DLW07}, Lemma 3.2.1. For each $n$ this can 
be made arbitrarily close to $1$ by choosing $a$ close to $1$.

\subsection{Generalized Pareto mixture}
Recall from Example \ref{exmp:gpd} that the GPD mixture algorithm has,
with $a \in (0,1)$, 
  \begin{align*}
    g_i(x \mid s) &= \alpha[a(b-s)]^\alpha x^{-\alpha-1} I\{x > a(b-s)\},
    \quad 1 \leq i \leq n-1.\\
    g_n(x \mid s) &= \alpha
    (b-s)^\alpha x^{-\alpha-1} I\{x > b-s\}I\{s \leq b-b(1-a)^{n-1}\}\\ 
    & \quad + f(x)I\{s > b-b(1-a)^{n-1}\}.
  \end{align*}
  First we check \eqref{eq:main} and \eqref{eq:sec}. 
  Karamata's theorem implies
  $\alpha \overline{F}(b) \sim b f(b)$. Then, for any $s < 1$,
  \begin{align}\label{eq:gpdcondratio}
    \frac{g_i(bx \mid bs)}{f(bx)} \overline{F}(b) &= 
    \frac{\alpha (bx)^{-\alpha-1} (ba(1-s))^\alpha \overline{F}(b)}{f(bx)}
    \nonumber \\
    &=  \frac{\alpha \overline{F}(bx)}{bx f(bx)} \frac{a^\alpha(1-s)^\alpha
        \overline{F}(b)}{x^\alpha \overline{F}(bx)} \to
    a^\alpha(1-s)^\alpha.
  \end{align}
  uniformly for $x \geq a(1-s)$. In particular \eqref{eq:main} is
  satisfied. Since
  \begin{align*}
    \frac{f(by)}{g_n(by\mid b s)} &= \frac{b f(bx)}{\alpha
    (1-s)^\alpha y^{-\alpha-1}} I\{y > 1-s\}I\{s \leq 1-(1-a)^{n-1}\} \\
    & \quad + I\{s > 1-(1-a)^{n-1}\}
  \end{align*}
  is bounded on $s \leq 1$, $y > 1-s$, \eqref{eq:sec} also holds. By
  Lemma \ref{lem:RNbddgen} the normalized Radon-Nikodym derivative is
  bounded. 
  By the arguments above
  \eqref{eq:thm} holds with $h_i(y \mid s) =
  a^{-\alpha}(1-s)^{-\alpha}$, $1 \leq i \leq n-1$ and $h_n(y \mid s)
  = (1-s)^{-\alpha}$. It follows by
  Theorem \ref{thm:main} that 
  \begin{align*}
    \lim_{b\to \infty} \frac{E \hat p_b^2}{p_b^2} = n^{-2}\Big(\prod_{i=1}^{n-1}
    \frac{1}{p_i} + a^{-\alpha} \sum_{j=1}^{n-1}
    \frac{1}{q_j}\prod_{i=1}^{j-1}\frac{1}{p_i}\Big).
  \end{align*}
  This is identical to \eqref{eq:condmix}, so $p_i$ can be chosen
  according to \eqref{eq:p} to minimize the relative error.
  
\subsection{Scaling mixtures}\label{sec:scalingmix}
In the scaling mixture algorithm presented in Example
\ref{exmp:scalemix} the large variables are generated by sampling from
the original density and multiplying
with a large number. In this section we study some variations of this
algorithm. Recall that, in the context of scaling mixtures,  always
assume that the orginal density $f$ is strictly positive on
$(0,\infty)$.

The first scaling mixture algorithm, called {\it scaling mixture I}, is
constructed as follows. Write $f(x) = x^{-\alpha-1}L(x)$ with $L$
slowly varying. Suppose $\inf_{x>x_0}L(x) =: L_* > 0$ for some $x_0>0$. 
The scaling mixture algorithm, with $\lambda > 0$, has
\begin{align*}
  g_i(x\mid s) &= (\lambda b)^{-1}f(x/\lambda b)I\{x > 0\} + f(x)I\{x
  \leq 0\},\quad i=1,\dots,n-1,\\
  g_n(x \mid s) &=  (\lambda b)^{-1}f(x/\lambda b)I\{x > 0, s \leq
  b-b(1-a)^{n-1}\}\\  
  & \quad + f(x)I\{x \leq 0 \text{ or } s > b-b(1-a)^{n-1}\}.
\end{align*}
To generate a sample $X$ from $g_i$ proceed as follows. Generate a candidate
$X'$ from $f$. If $X' \leq 0$ put $X = X'$ and if $X' > 0$, put $X =
\lambda b X'$. 

Take $a\in(0,1)$, using Karamata's theorem, $\alpha \overline{F}(b) \sim b f(b)$, 
we have, for $1 \leq i \leq n$, and $s \leq 1-(1-a)^{i-1}$,
\begin{align*}
  \frac{g_i(bx \mid bs)}{f(bx)}\overline{F}(b) &= 
  \frac{f(\frac{x}{\lambda})}{\lambda b f(bx)}\overline{F}(b)I\{x >
  0\} + \overline{F}(b)I\{x \leq 0\}\\
  &= \frac{x f(\frac{x}{\lambda})}{\alpha \lambda} \frac{\alpha
      \overline{F}(bx)}{bx
      f(bx)}\frac{\overline{F}(b)}{\overline{F}(bx)}I\{x >   
  0\} + \overline{F}(b)I\{x \leq 0\}\\
 & \to
   \frac{ x^{\alpha+1} f(\frac{x}{\lambda})}{\alpha
    \lambda}  
\end{align*}
uniformly for $x \geq 1-s$. Since $x^{\alpha+1}f(x/\lambda) \geq
  \lambda^{\alpha+1}L_*>0$, the condition \eqref{eq:main} holds. Note,
  however, that \eqref{eq:main} fails if $L_* = 0$. 
Since
\begin{align*}
  \frac{f(bx)}{g_n(bx \mid bs)} = \frac{\lambda b
    f(bx)}{f(\frac{x}{\lambda})}I\{x > 0, s \leq 1-(1-a)^{n-1}\} +
I\{x \leq 0 + s > 1-(1-a)^{n-1}\} 
\end{align*}
is bounded on $s \leq 1$, $x > 1-s$ condition \eqref{eq:sec} also
holds. From the calculation above we see that
  \eqref{eq:thm} is satisfied
  with $h(x \mid s) = \alpha \lambda [x^{\alpha+1}f(x/\lambda)]^{-1}$.  
  In particular, the asymptotic upper bound for the normalized second
  moment is
  \begin{align*}
   \frac{1}{n^2} \int R(y)dm =  \frac{1}{n^2} \lambda^{-2\alpha}\int_{1/\lambda}^\infty
    \frac{\alpha^2}{x^{2(\alpha+1)}f(x)} dx \sum_{i=1}^{n} \frac{1}{q_i} 
    \prod_{j=1}^{i-1}\frac{1}{p_j},
  \end{align*}
  with $q_n = 1$. It is straightforward to check that $\frac{1}{n^2}\sum_{i=1}^n \frac{1}{q_i} 
  \prod_{j=1}^{i-1}\frac{1}{p_j}$ is minimized at 
  \begin{align*}
    p_i = 1-\frac{1}{n- i +1},\quad q_i
    = 1-p_i,
  \end{align*}
  with minimum equal to $1$. The parameter $\lambda$ can be
  chosen to control the factor
  \begin{align*}
    \lambda^{-2\alpha}\int_{1/\lambda}^\infty
    \frac{\alpha^2}{x^{2(\alpha+1)}f(x)} dx. 
  \end{align*}
In some cases this can be minimized analytically. 
\begin{exmp}
  Consider a Pareto density of the form $f(x) = \alpha(1+x)^{-\alpha -
    1}$, $x > 0$.  Then 
\begin{align*}
  \lambda^{-2\alpha}\int_{1/\lambda}^\infty
  \frac{\alpha^2}{x^{2(\alpha+1)}f(x)} dx =
    \lambda^{-2\alpha}\int_{1/\lambda}^\infty
    \alpha\Big(\frac{1+x}{x^2}\Big)^{\alpha+1} dx. 
  \end{align*}
  If $\alpha = 1$ this is minimized at $\lambda = \sqrt{3}$ with
  minimum $\frac{2+\sqrt{3}}{\sqrt{3}}$. 
\end{exmp}

In the scaling mixture algorithm we assume $L_* > 0$. This rules
out distributions whose slowly varying function tends to $0$. However,
this is not a severe problem. One way to avoid it is to slightly modify
the previous algorithm. The {\it scaling mixture II} algorithm has, with
$\lambda > 0$, $u \in (0,1)$, $\delta > 0$, and $a \in (0,1)$,
\begin{align*}
  g_i(x\mid s)  &= g(x) \\ & = (\lambda b)^{-1}f(x/\lambda
  b)I\{0 < x \leq \lambda b u\} 
    \\ & \quad +
  \frac{1}{(1+\delta)\lambda b}\Big(\frac{x}{\lambda
    b}\Big)^{\frac{1}{1+\delta}-1} f([x/\lambda
  b]^{\frac{1}{1+\delta}})I\{x \geq \lambda b u^{1+\delta}\} +
  f(x)I\{x \leq 0\},\\
g_n(x\mid s) &= g(x)I\{s \leq b-b(1-a)^{n-1}\} + f(x)I\{s >
b-b(1-a)^{n-1}\}. 
\end{align*}
The density $g_i$ is based on the following sampling procedure. 
To generate a sample $X$ from $g_i$, first generate a candidate
$X'$ from $f$. If $X' \leq 0$ put $X = X'$, if $0 < X' \leq u$, put $X =
\lambda b X'$, and if $X' > u$ put $X = \lambda b (X')^{1+\delta}$. 

Similar to the scaling mixture I algorithm it follows that, for $1 \leq i \leq n$,
\begin{align*}
  &\frac{g_i(bx \mid bs)}{f(bx)}\overline{F}(b) \\ &
  \quad \to
    x^{\alpha+1} \frac{f(\frac{x}{\lambda})}{\alpha
    \lambda}I\{0< x \leq \lambda u\} +  
    \frac{ x^{\frac{1}{1+\delta}+\alpha}  
      f([x/\lambda]^{\frac{1}{1+\delta}})}{(1+\delta)\alpha
      \lambda^{\frac{1}{1+\delta}}}I\{x \geq \lambda u^{1+\delta}\}    
  \end{align*}
  uniformly for $x \geq 1-s$ and $s$. Since 
  \begin{align*}
    x^{\frac{1}{1+\delta}+\alpha}f([x/\lambda]^{\frac{1}{1+\delta}}) = 
    \lambda^{-\frac{\alpha+1}{1+\delta}} x^{\alpha(1-\frac{1}{1+\delta})}L(x/\lambda),
  \end{align*}
  is bounded from below for $x \geq 1-s$ \eqref{eq:main} holds. 
  Just as for the scaling mixture I algorithm \eqref{eq:sec} also
  holds. \eqref{eq:thm} hold with 
  \begin{align*}
    h_i(y\mid s)=\frac{\alpha^2\lambda}{y^{2\alpha+2}f(y/\lambda)}+\frac{(1+\delta)\alpha^2
      \lambda^{\frac{1}{1+\delta}}}{y^{2\alpha+\frac{1}{1+\delta}+1}f([y/\lambda]^{\frac{1}{1+\delta}})}.
  \end{align*}
  The asymptotic upper bound for the normalized second moment is hence
  \begin{align}\label{eq:sm2re}
    \int R(y)dm &=  \sum_{i=1}^{n} \frac{1}{q_i} 
    \prod_{j=1}^{i-1}\frac{1}{p_j}\nonumber \\
    & \quad \times \Big(\int_1^{\lambda
      u}\frac{\alpha^2\lambda}{x^{2\alpha+2}f(x/\lambda)} dx + 
    \int_{\lambda u^{1+\delta}}^\infty \frac{(1+\delta)\alpha^2
      \lambda^{\frac{1}{1+\delta}}}{x^{2\alpha+\frac{1}{1+\delta}+1}f([x/\lambda]^{\frac{1}{1+\delta}})} dx\Big). 
  \end{align}
  with $q_n = 1$. As above $\frac{1}{n^2}\sum_{i=1}^n \frac{1}{q_i} 
  \prod_{j=1}^{i-1}\frac{1}{p_j}$ is minimized at 
  \begin{align*}
    p_i = 1-\frac{1}{n- i +1},\quad q_i
    = 1-p_i,
  \end{align*}
  with minimum equal to $1$. The remaining parameters $\lambda$ and
  $u$ can be chosen to control the integrals in \eqref{eq:sm2re}.

\section{Achieving asymptotically optimal relative error} \label{sec:asopt}

In the previous section we observed that the
conditional mixture algorithm and the GPD mixture algorithm can be
designed to have almost asymptotically optimal relative error. A small
asymptotic relative error is obtained by choosing the parameter $a$
close to $1$. In this section we prove that these algorithms have 
asymptotically optimal relative error. This is accomplished by letting the parameter
$a$ depend on the threshold $b$ in such a way that
$a\to 1$ slowly as $b\to \infty$. For simplicity, we assume that $X_1>0$ throughout this
section.
\begin{thm}\label{thm:asopt}
  Let $\nu_n^b$ be the measure defined by the conditional mixture
  algorithm. Let $p_i = \frac{n-i}{n-i+1}$, $q_i
  = 1-p_i$, and assume that $1-a=1-a_b\sim
  \mathcal{O} (b^{-\frac{1}{2(n-1)}+\delta})$ for some $0<\delta< \frac{1}{2(n-1)}$.  
  Then, the conditional mixture algorithm has asymptotically
  optimal relative error for computing $p_b$. That is, 
  \begin{align*}
    \lim_{b\to\infty} \frac{E\hat p_b^2}{p_b^2}=1.
  \end{align*}
\end{thm}
\begin{rem}
  In Theorem \ref{thm:asopt} the conditional
  mixture algorithm can be replaced by the GPD mixture algorithm.
\end{rem}
\begin{proof}
  First rewrite the normalized second moment as in \eqref{eq:nsm}:
  \begin{align*}
    \frac{E \hat p_b^2}{p_b^2} =  
    \frac{1}{n^2}\int\limits_{s_n>1}
    \frac{1}{\overline{F}(b)}\frac{d\mu_n}{d\nu_n^b}(by)m_b(dy) =
    \frac{1}{n^2}\int\limits_{s_n > 1} 
    R_b(y) m_b(dy).   
  \end{align*}
  Fix $a_0 \in (0,1)$. Define the sets
  \begin{align*}
    A_i &= \{y \in\R^n : y_i > 1-s_{i-1}, y_j \leq 1-s_{j-1}, j < i\},\\
    B_i &= \{y \in\R^n : y_i \leq a_0(1-s_{i-1})\},\\
    C_i &= \{y \in\R^n : a_0(1-s_{i-1}) < y_i \leq a(1-s_{i-1})\},\\
    D_i &= \{y \in\R^n :  a(1-s_{i-1}) < y_i \leq 1-s_{i-1}\}.
  \end{align*}
  Then $\{s_n > 1\} \subset \cup_{i=1}^n A_i$ and each $A_i$ can be
  written as the disjoint union of the $3^{i-1}$ sets of the form
  \begin{align}\label{eq:int}
    I_1 \cap I_2 \cap \dots \cap I_{i-1}\cap A_i,
  \end{align}
  where each $I_j$ is either $B_j$, $C_j$, or $D_j$.  Each
  intersection \eqref{eq:int} is of one of the types below. 
  \begin{itemize}
  \item[(i)] $I_j = B_j$ for each $j=1,\dots, i-1$. 
  \item[(ii)] among the sets $I_1, \dots, I_{i-1}$ there is
    at least one $j$ for which $I_j = C_j$ and no $j$ with $I_j =
    D_j$.
  \item[(iii)] among the sets $I_1, \dots, I_{i-1}$ there is
    at least one $j$ for which $I_j = D_j$. 
  \end{itemize}
Next we treat the integrals  
\begin{align*}
  \frac{1}{n^2}\int\limits_{I_1\cap \dots \cap I_{i-1}\cap A_i} 
  R_b(y) m_b(dy),
\end{align*}
separately. The intersection belongs to one of the three types. 

Type (i): Consider $y \in B_1 \cap \dots \cap B_{i-1}\cap A_i$. Then
$s_{i-1} \leq 1-(1-a_0)^{i-1}$ and 
\begin{align*}
  R_b(y) &\leq \prod_{j=1}^{i-1}\frac{1}{p_j} \times \frac{1}{p_i
    \overline{F}(b) + q_i
    \frac{\overline{F}(b)}{\overline{F}(ba(1-s_{i-1}))}} \leq \prod_{j=1}^{i-1}\frac{1}{p_j}
 \frac{1}{q_i}\frac{\overline{F}(ba(1-s_{i-1}))}{\overline{F}(b)}. 
\end{align*}
Fix arbitrary $\vep > 0$. Then, for $b$ sufficiently large, $a_b >
1-\vep$ and the expression in the last display is bounded above by
\begin{align*}
\prod_{j=1}^{i-1}\frac{1}{p_j}
 \frac{1}{q_i}\frac{\overline{F}(b(1-\vep)(1-s_{i-1}))}{\overline{F}(b)}=:
 R^*_b(y).  
\end{align*}
It follows that $R^*_b(y) \to \prod_{j=1}^{i-1}\frac{1}{p_j}
 \frac{1}{q_i}(1-\vep)^{-\alpha}(1-s_{i-1})^{-\alpha}$ uniformly on
 $B_1 \cap \dots \cap B_{i-1} \cap A_i$ and then it follows by the
 arguments in the proof of Theorem \ref{thm:main} that 
 \begin{align*}
    \limsup_{b \to \infty} \frac{1}{n^2}\int\limits_{B_1\cap \dots \cap B_{i-1}\cap A_i} 
  R_b(y) m_b(dy) \leq \frac{1}{n^2} \prod_{j=1}^{i-1}\frac{1}{p_j}
 \frac{1}{q_i} (1-\vep)^{-\alpha}. 
 \end{align*}
Since $\vep > 0$ was arbitrary we can let $\vep \to 0$ to get
 \begin{align*}
    \limsup_{b \to \infty} \frac{1}{n^2}\int\limits_{B_1\cap \dots \cap B_{i-1}\cap A_i} 
  R_b(y) m_b(dy) \leq \frac{1}{n^2} \prod_{j=1}^{i-1}\frac{1}{p_j}
 \frac{1}{q_i}. 
 \end{align*}
Type (ii): For $y \in I_1 \cap I_2 \cap \dots \cap I_{i-1}\cap A_i$ it
holds that $s_{i-1} \leq 1-(1-a)^{i-1}$. Proceeding as in the Type (i)
case, for $\vep > 0$ and $b$ sufficiently large,  
\begin{align*}
  R_b(y) &\leq \prod_{j=1}^{i-1}\frac{1}{p_j} \times \frac{1}{p_i
    \overline{F}(b) + q_i
    \frac{\overline{F}(b)}{\overline{F}(ba(1-s_{i-1}))}} \\ & \leq
  \prod_{j=1}^{i-1}\frac{1}{p_j} 
  \frac{1}{q_i}\frac{\overline{F}(ba(1-s_{i-1}))}{\overline{F}(b)}\\ 
  & \leq
  \prod_{j=1}^{i-1}\frac{1}{p_j} 
  \frac{1}{q_i}\frac{\overline{F}(b(1-\vep)(1-a)^{i-1})}{\overline{F}(b)}.  
\end{align*}
It follows that
\begin{align}
  \frac{1}{n^2}&\int\limits_{I_1\cap \dots \cap I_{i-1}\cap A_i} 
  R_b(y) m_b(dy) \nonumber \\ & \leq \frac{1}{n^2} \prod_{j=1}^{i-1}\frac{1}{p_j}
  \frac{1}{q_i}\frac{\overline{F}(b(1-\vep)(1-a)^{i-1})}{\overline{F}(b)}
  m_b(I_1\cap \dots \cap I_{i-1}\cap A_i). \label{eq:type2.1}
\end{align}
Let $k$ be the first index between $1$ and $i-1$ such that $I_k =
C_k$. Then $I_l = B_l$ for each $1 \leq l \leq k-1$, and $s_{k-1}\leq 1-(1-a_0)^{k-1}$, 
whereas $s_{i-1}\leq 1-(1-a)^{i-1}$, 
\begin{align}
  m_b(I_1\cap \dots \cap I_{i-1}\cap A_i) & \leq \frac{P(Y_k>ba_0(1-s_{k-1}), Y_i>b(1-s_{i-1}))}{\overline{F}(b)}\nonumber\\  
  & \leq \frac{P(Y_k>ba_0(1-a_0)^{k-1}, Y_i>b(1-a)^{i-1})}{\overline{F}(b)}\nonumber\\  
  &=\frac{\overline{F}(ba_0(1-a_0)^{k-1})
    \overline{F}(b (1-a)^{i-1})}{\overline{F}(b)}\nonumber \\ &\leq
  \frac{\overline{F}(ba_0(1-a_0)^{i-1}) 
    \overline{F}(b (1-a)^{i-1})}{\overline{F}(b)}\label{eq:opt}.
\end{align}
Putting this into \eqref{eq:type2.1} yields the upper bound
\begin{align*}
  \frac{1}{n^2}&\int\limits_{I_1\cap \dots \cap I_{i-1}\cap A_i} 
  R_b(y) m_b(dy) \\
  & \leq  \frac{1}{n^2} \prod_{j=1}^{i-1}\frac{1}{p_j}
   \frac{1}{q_i}\Big[\frac{\overline{F}(b(1-\vep)(1-a)^{i-1})}{\overline{F}(b)}\Big]^2\overline{F}(ba_0(1-a_0)^{i-1})
   \\
   & \leq  \frac{1}{n^2} \prod_{j=1}^{i-1}\frac{1}{p_j}
   \frac{1}{q_i}\Big[\frac{\overline{F}(b(1-\vep)(1-a)^{n-1})}{\overline{F}(b)}\Big]^2\overline{F}(ba_0(1-a_0)^{n-1}).
 \end{align*}
This converges to $0$ as $b \to \infty$ by the choice of $a = a_b$. 

Type (iii):  For $I_1 \cap I_2 \cap \dots \cap I_{i-1}\cap A_i$ of
type (iii) we let $j$ denote the first index for which $I_j =
D_j$. Suppose first that $I_k = B_k$ for each $k=1,\dots,j-1$. Then, 
$s_{j-1} \leq 1-(1-a_0)^{j-1}$ and, for arbitrary $\vep > 0$ and 
$b$ sufficiently large, 
\begin{align}\label{eq:type3.1}
  R_b(y) &\leq \prod_{k=1}^{j-1}\frac{1}{p_k}
 \frac{1}{q_j}\frac{\overline{F}(b(1-\vep)(1-a_0)^{j-1})}{\overline{F}(b)}, 
\end{align}
which is bounded in $b$. In addition,
\begin{align}
  & m_b(B_1\cap \dots \cap B_{j-1}\cap D_j)  \leq \frac{P(Y_j\in D_j,
    S_{j-1}\leq 1-(1-a_0)^{j-1})}{\overline{F}(b)}\nonumber\\ 
  &\quad\quad \leq \int\limits_{B_1\cap
    \dots \cap B_{j-1}} 
  \frac{\overline{F}(ba(1-s_{j-1}))-
    \overline{F}(b(1-s_{j-1}))}{\overline{F}(b)} \mu_n(dy) \to 0,\label{eq:type3.2}
\end{align}
as $b \to \infty$, by the bounded convergence theorem. Combining
\eqref{eq:type3.1} and \eqref{eq:type3.2} we see that
\begin{align*}
  \limsup_{b\to \infty} \int\limits_{B_1\cap
    \dots \cap B_{j-1}} R_b dm_b = 0.
\end{align*}
Finally, suppose $I_k = C_k$ for some $k = 1,\dots,j-1$. 
Then, $s_{j-1} \leq 1-(1-a)^{j-1}$ and, for arbitrary $\vep > 0$ and 
$b$ sufficiently large, 
\begin{align}\label{eq:type3.3}
  R_b(y) &\leq \prod_{k=1}^{j-1}\frac{1}{p_k}
 \frac{1}{q_j}\frac{\overline{F}(b(1-\vep)(1-a)^{j-1})}{\overline{F}(b)}.
\end{align}
In addition, just as in \eqref{eq:opt},
\begin{align}\label{eq:type3.4}
  m_b(I_1\cap \dots \cap I_{j-1}\cap D_j) & \leq 
  \frac{\overline{F}(ba_0(1-a_0)^{j-1}) 
    \overline{F}(b(1-\vep)(1-a)^{j-1})}{\overline{F}(b)}.
\end{align}
Combining \eqref{eq:type3.3} and \eqref{eq:type3.4} we see that
\begin{align*}
  & \limsup_{b\to \infty} \int\limits_{B_1\cap
    \dots \cap B_{j-1}} R_b dm_b \\ & \quad \leq \limsup_{b \to \infty}
  \prod_{k=1}^{j-1}\frac{1}{p_k} 
 \frac{1}{q_j}\Big[\frac{\overline{F}(b(1-\vep)(1-a)^{j-1})}{\overline{F}(b)}\Big]^2\overline{F}(ba_0(1-a_0)^{j-1})
 = 0, 
\end{align*}
by the choice of $a = a_b$. 
\end{proof}

\section{Numerical illustrations}

In this section we examine the performance of the scaling mixture algorithm, referred to as the SM algorithm.
We perform a preliminary  test using Pareto-distributed positive random variables and 
compare the algorithm with the conditional mixture algorithm in
\cite{DLW07}, which we refer to as the DLW algorithm, and the conditional Monte Carlo algorithm
in \cite{AK06}. For comparision, we first consider the same setting as 
in \cite{DLW07}, Table IV, pp. 18. The so-called true value in Table \ref{table1} 
was obtained from the same table.
Each estimate was calculated using $N=10^4$ samples of $S_n$. 
This estimation was repeated 100 times and the mean estimate, the mean standard
error and the mean calculation time were calculated. 
The parameter $a$ in the DLW algorithm was chosen equal to 0.999 and the
parameter $\lambda$ in the scaling mixture algorithm was chosen equal to 1 in
Table \ref{table1} and equal to $\sqrt{3}$, the optimal value, in Table \ref{table2}.

The standard Monte Carlo estimation is inferior to both importance sampling algorithms.
The conditional Monte Carlo algorithm performs best for most probabilites in this study. 

\begin{table}[ht]
\caption{Simulations of $P(S_{n}>b)$, 
where $S_n=\sum_{i=1}^n X_i$ and $P(X_1>x)=(1+x)^{-1/2}$, $a=0.999$ and $\lambda=1$.
$N=10^4$ samples were used for each estimation, repeated 100 times.. 
}\label{table1}
 {\scriptsize
\begin{tabular}{|c|c|c|c|c|c|c|c|}
\hline
  $n$ 	&  $b$ 		&   True value	&  MIS   	&  DLW		&   CMC		&   MC 		&\\
\hline
    5	&  5e+05 	&   0.007071	&   0.0070744	&  0.0070714	&  0.00707034	&  0.0069960	&  Avg.~est.\\
	&   		&  	 	&   (7.26e-05)	& (6.10e-06) 	&  (4.89e-06)	& (4.88e-05)	& (A.~std.~err.)\\
	&		&		&   [0.816]	&  [0.799]	&   [0.731]	&  [0.685]	& [A.~time (s)]\\
\hline 
	&  5e+11	&   7.0711e-06	&   7.0776e-06	&  7.0710e-06	&  7.0711e-06	&  1.8000e-05	&\\
	&  		&   	 	&  (7.53e-08)	& (1.86e-09)	&  (2.71e-11)	&  (1.56e-05)	&\\
	&		&		&  [1.005]	&   [0.990]	&   [0.908]	&   [0.840]	&\\
\hline
  15 	&  5e+05	&   0.02121	&   0.021188 	&   0.021215 	&   0.021210 	&   0.021724 	&\\
	&   		&   		&  (2.07e-04)	& (4.15e-05)	&  (2.72e-05)	&  (2.05e-03)	&\\
	&		&		&   [1.224]	&  [1.219 ]	&    [1.092]	&   [1.006]	&\\
\hline
	&  5e+11	&   2.1213e-05	&   2.1224e-05 	&   2.1214e-05	&  2.1213e-05	&  1.800e-05	&\\
	&  		&  		&   (2.25e-07)	&  (5.82e-09)	&  (3.09e-10)	&  (1.80e-05) 	&\\
	&		&		&    [1.450]	&    [1.456]	&   [1.283] 	&   [1.179]	&\\
\hline
  25 	&  5e+05	&   0.035339	&   0.035330	&  0.035348 	&  0.035347	&   0.035462	&\\
	&   		&    		&   (3.32e-04)	&  (9.06e-05)	&  (5.89e-05)	&  (2.61e-03)	&\\
	&		&		&   [1.712]	& [1.729]	&   [1.478]	&   [1.366]	&\\
\hline
	&  5e+11	&   3.5355e-05	&   3.5338e-05	&  3.5355e-05	&  3.5355e-05 	&  3.8000e-05	&\\
	&  		&   		&  (3.77e-07)	&  (1.04e-09)	&  (1.32e-09)	&  (3.68e-05) 	&\\
	&		&		&    [1.993]	&   [2.016]	&   [1.689]	&   [1.559]	&\\
\hline

\end{tabular} }
\end{table}

\begin{table}[ht]
\caption{Simulations of $P(S_{n}>b)$, 
where $S_n=\sum_{i=1}^n X_i$, $P(X_1>x)=(1+x)^{-1}$, $a=0.999$ and $\lambda=\sqrt{3}$.
$N=10^4$ samples were used for each estimation, repeated 100 times.
}\label{table2}
  {\scriptsize
\begin{tabular}{|c|c|c|c|c|c|c|c|}
\hline
  $n$ 	&  $b$ 		&   True value	&  MIS   	&  DLW		&   CMC		&   MC 		&\\
\hline
    5	&  5e+05 	&   1.0001e-05	&   1.0020e-05	&  1.0001e-05	&  1.0001e-05	&  6.000e-06	&  Avg.~est.\\
	&   		&  	 	&   (1.07e-07)	& (2.78e-09) 	&  (2.58e-10)	& (6.00e-6)	& (std.~err.)\\
	&		&		&   [0.429]	&  [0.415]	&   [0.433]	&  [0.346]	& [time (s)]\\
\hline 
	&  5e+11	&   1.0000e-13	&   9.9996e-12	&  9.9999e-12	&  1.0000e-13	&    0		&\\
	&  		&   	 	&  (1.07e-13)	& (2.79e-15)	&  (8.59e-22)	&    (0)	&\\
	&		&		&  [0.433]	&   [0.418]	&   [0.430]	&   [0.352]	&\\
\hline
  15 	&  5e+05	&   3.0010e-05	&   3.0004e-05	&  3.0011e-05 	&   3.0010e-05	&   3.0000e-05 	&\\
	&   		&   		&  (3.21e-07)	& (1.12e-08)	&  (1.74e-09)	&  (2.71e-05)	&\\
	&		&		&   [0.491]	&  [0.445]	&    [0.437]	&   [0.375]	&\\
\hline
	&  5e+11	&   3.0000e-11	&   2.9990e-11 	&   3.0000e-11	&  3.0000e-11	&  	0	&\\
	&  		&  		&   (3.22e-13)	&  (9.06e-15)	&  (1.75e-20)	&  	(0) 	&\\
	&		&		&    [0.490]	&    [0.445]	&   [0.431] 	&   [0.365]	&\\
\hline
  25 	&  5e+05	&  5.0029e-05	&   5.0098e-05	&  5.00274e-05 	&  5.00290e-05	&   3.7000e-05	&\\
	&   		&    		&   (5.37e-07)	&  (1.90e-08)	&  (4.10e-09)	&  (3.34e-05)	&\\
	&		&		&   [0.561]	& [0.485]	&   [0.432]	&   [0.386]	&\\
\hline
	&  5e+11	&   5.0000e-11	&   4.9970e-11	&  4.9998e-11	&  5.0000e-11 	&  0		&\\
	&  		&   		&  (5.38e-13)	&  (1.65e-14)	&  (1.54e-20)	&  (0) 		&\\
	&		&		&    [0.556]	&   [0.479]	&   [0.439]	&   [0.382]	&\\
\hline

\end{tabular} }
\end{table}

\end{document}